\documentclass[11pt]{article}

\usepackage[utf8]{inputenc}
\usepackage[english]{babel}
\usepackage{amsmath, amsthm, amssymb, amscd}

\usepackage{graphicx}
\usepackage{authblk}

\newcommand{\Bin}{\operatorname{Bin}}
\newcommand{\Beta}{\operatorname{Beta}}
\newcommand{\data}{\mathbf{X}^k}
\newtheorem{theorem}{Theorem}

\begin{document}

\title{Posterior Consistency in the Binomial $(n,p)$ Model with Unknown $n$ and $p$: A Numerical Study}

\author[1]{Laura Fee Schneider\thanks{laura-fee.schneider@mathematik.uni-goettingen.de}}
\author[1]{Thomas Staudt\thanks{	thomas.staudt@stud.uni-goettingen.de}}
\author[1,3]{Axel Munk\thanks{munk@math.uni-goettingen.de}}
\affil[1]{Institute for Mathematical Stochastics, University of G\"{o}ttingen}
\affil[3]{Max Planck Institute for Biophysical Chemistry, G\"{o}ttingen}

\date{ }
\maketitle

\abstract{ Estimating the parameters from $k$ independent Bin$(n,p)$ random variables, when both parameters $n$ and $p$ are unknown, is relevant to a variety of applications. It is particularly difficult if $n$ is large and $p$ is small. Over the past decades, several articles have proposed Bayesian approaches to estimate $n$ in this setting, but asymptotic results could only be established recently in \cite{Schneider}. There, posterior contraction for $n$ is proven in the problematic parameter regime where $n\rightarrow\infty$ and $p\rightarrow0$ at certain rates. In this article, we study numerically how far the theoretical upper bound on $n$ can be relaxed in simulations without losing posterior consistency.}

\section{Introduction}
\label{sec:1}
We consider estimating the parameter $n$ of the binomial distribution from $k$ independent observations when the success probability $p$ is unknown. This situation is relevant in many applications, for example in estimating the population size of a species \cite{Raftery} or the total number of defective appliances \cite{Draper-Guttman}.
Another recent application is quantitative nanoscopy, see \cite{Schneider}. There, the total number of fluorescent markers (fluorophores) attached to so-called DNA-origami is estimated from a time series of microscopic images. The number of active fluorophores counted in each image is modeled as binomial observation, where the probability $p$ that a fluorophore is active in the respective image is very small (often below $5\%$). 

This setting, where the success probability $p$ is small (and $n$ potentially large), is very challenging. The difficulties that arise can be understood by considering the following property of the binomial distribution: if $n$ converges to infinity, $p$ converges to zero, and the product $np$ converges to $\lambda > 0$, then a Bin$(n,p)$ random variable converges in distribution to a Poisson variable with parameter $\lambda$. Thus, the binomial distribution converges to a distribution with a single parameter. This suggests that it gets harder to derive information about the two parameters separately when $n$ is large and $p$ small.

In this context, it is instructive to look at the sample maximum $M_k$ as an estimator for $n$, which was suggested by Fisher in 1941 \cite{Fisher}. Although it turns out to be impractical, see \cite{DasGupta}, the sample maximum is consistent and converges in probability for fixed parameters $(n, p)$ exponentially fast to the true $n$, as $k\rightarrow\infty$. This can be seen from
\begin{equation}\label{prob-max}
  \mathbf{P}\left( M_k =n  \right)= 1-(1-p^n)^k,
\end{equation}
which implies, by Bernoulli inequality and since $1-x\leq e^{-x}$, that
\begin{equation*}
1-e^{-kp^n} \leq \mathbf{P}\left( M_k =n  \right)\leq kp^n.
\end{equation*}
In an asymptotic setting where $n\rightarrow\infty$ and $p\rightarrow0$ such that $kp^n\rightarrow0$, the probability in \eqref{prob-max} no longer converges to one. Thus, the sample maximum is a consistent estimator for $n$ only as long as $kp^n\rightarrow\infty$. The condition $e^n=O(k)$ is necessary for this to hold.\\

Estimating $n$ in this difficult regime becomes more manageable by including prior knowledge about $p$. We therefore consider random $N$ and $P$, and variables $X_1,\dots, X_k$ that are independently $\Bin(n,p)$ distributed given that $N=n$ and $P=p$.
Various Bayesian estimators have been suggested over the last 50 years, see \cite{Draper-Guttman,Raftery,Berger,Chilko,Hamedani}. In all of this work, a product prior for $(N,P)$ is used, and the prior $\Pi_P$ on $P$ is chosen as beta distribution Beta$(a,b)$ for some $a,b > 0$. 
Since this is the conjugate prior, it is a natural choice. In contrast, there is quite some discussion about the most suitable prior $\Pi_N$ for $N$, see for example \cite{Kahn,Link,Villa,Berger}. Therefore, the asymptotic results in \cite{Schneider} are described flexible in terms of $\Pi_N$, and they only require a condition that ensures that enough weight is put on large values of $n$ (see equation \eqref{tail-bound} in Section \ref{sec:2}).

In \cite{Schneider}, we also introduce a new class of Bayesian point estimators for $n$, which we call scale estimators. 
We choose $\Pi_P\sim\Beta(a,b)$ and set $\Pi_N(m)\propto m^{-\gamma}$ for a positive value $\gamma$. If $\gamma>1$, the prior $\Pi_N$ is a proper probability distribution, but it is sufficient to ensure $\gamma +a>1$ in order to obtain a well-defined posterior distribution.
The scale estimator is then defined as the minimizer of the Bayes risk with respect to the relative quadratic loss,
$ l(x,y)=( x/y-1)^2. $
Following \cite{Raftery}, it is given by
\begin{equation}\label{E:ScaleBayes}
\hat{n}:=\frac{\mathbb{E} \left[\frac{1}{N}|\data\right]}{\mathbb{E} \left[\frac{1}{N^2}|\data\right]}=\frac{\sum_{m=M_k}^\infty \frac{1}{m}L_{a,b}(m)\Pi_N(m)}{\sum_{m=M_k}^\infty \frac{1}{m^2}L_{a,b}(m)\Pi_N(m)},
\end{equation}
where $\mathbf{X}^k=(X_1,\dots,X_k)$ denotes the sample, $M_k$ is the sample maximum, and $L_{a,b}$ is the beta-binomial likelihood, see \cite{Carroll-Lombard}. We refer to \cite{Schneider} for a detailed discussion and numerical study of this estimator.\\

The present article is structured as follows. In Section \ref{sec:2}, the main theorem (proven in \cite{Schneider}) is presented, which shows uniform posterior contraction in the introduced Bayes setting for suitable asymptotics of $n$ and $p$. 
The theorem states that $n^{6+\epsilon}=O(k)$ for $\epsilon>0$ is already sufficient for consistency of the Bayes estimator, improving significantly over the sample maximum.
In Section \ref{sec:3}, we then conduct a simulation study to closer investigate the restrictions for the parameters $n$ and $p$ needed to ensure consistency. Our findings indicate that estimation of $n$ is still consistent if $n^5=O(k)$, but that it becomes inconsistent for $n^3=O(k)$. 
It is hard to pin down the exact transition from consistency to inconsistency when $n^{\alpha}=O(k)$, but our results suggest that it happens close to $\alpha = 4$.
We discuss our results and provide several remarks in Section \ref{sec:4}.

\section{Posterior Contraction for $n$}
\label{sec:2}

To study posterior contraction in the binomial model we consider the Bayesian setting described in Section \ref{sec:1}. For fixed parameters $n$ and $p$ that are independent of the number of observations $k$, posterior consistency follows from Doob's theorem, see, e.g., \cite{vanderVaart-Asymptotic}. We extend this result to the class of parameters
\begin{equation} \label{M}
  \mathcal{M}_{\lambda}:=\left\{ (n_k,p_k)_k: 1/\lambda \leq n_kp_k\leq \lambda,\ n_k\leq \lambda\sqrt[6]{k/\log(k)} \right\}
\end{equation}
for fixed $\lambda > 1$.
Since we want to handle a variety of suitable prior distributions for $N$, we only require that $\Pi_N$ is a proper probability distribution on $\mathbb{N}$ that fulfills the condition
\begin{equation}\label{tail-bound}
  \Pi_N(m) \geq \beta e^{-\alpha m^2} \quad \forall\, m\in\mathbb{N}
\end{equation}
for some positive constants $\alpha$ and $\beta$.

\begin{theorem}[see \cite{Schneider}]\label{main_result}
  Conditionally on $N=n_k$ and $P=p_k$, let $X_1,\dots,X_k\overset{i.i.d.}{\sim} \Bin(n_k,p_k)$. For any prior distribution $\Pi_{(N,P)} = \Pi_N \Pi_P$ on $(N,P)$ with $\Pi_P=\Beta(a,b)$ for $a,b>0$, and where $\Pi_N$ satisfies \eqref{tail-bound}, we have uniform posterior contraction over the set $\mathcal{M}_{\lambda}$ of sequences $(n_k, p_k)_k$ defined in \eqref{M} for any $\lambda>1$, i.e.,
\begin{equation*}\label{E:Result_sup}
\sup_{(n_k, p_k)_k \in \mathcal{M}_{\lambda}} \mathbb{E}_{n_k,p_k}\big[ \Pi\big( N \neq n_k\ |\, \data \big) \big] \rightarrow 0,\ \text{as}~ k\rightarrow\infty.
\end{equation*}
\end{theorem}

This result directly implies consistency of the scale estimator (\ref{E:ScaleBayes}) for parameter sequences in $\mathcal{M}_{\lambda}$. The flexible restrictions on the prior distribution allow to apply the result to the estimators derived in \cite{Chilko} and \cite{Hamedani} as well. Furthermore, it is possible to extend the statement of Theorem \ref{main_result} to improper priors on $N$, as done in Theorem 2 in \cite{Schneider}, in order to cover the estimators in \cite{Draper-Guttman} and \cite{Berger}.

\section{Simulation Study}
\label{sec:3}

The theorem presented in the previous section states that the asymptotic behavior $n_k \sim O\big(\sqrt[6]{k/\log(k)}\big)$ leads to posterior contraction of $N$ for suitable priors, as long as $n_k p_k$ stays in a compact interval bounded away from zero. In this section we try to answer the question by how much the constraints on $\mathcal{M}_\lambda$ in Theorem \ref{main_result} can be relaxed.
We address this problem by studying the relation between posterior contraction and the order $\alpha > 0$ when $n_k \sim O\big(\sqrt[\alpha]{k}\big)$. More precisely, we are interested in the smallest $\alpha = \alpha^*$ such that the result 
\begin{equation}\label{eq:contraction}
  \mathbb{E}_{n_k, p_k}\left[ \Pi\left( N\neq n_k\,|\, \mathbf{X}^k \right)  \right] \rightarrow 0,\ ~\text{as}~ \ k\rightarrow\infty,
\end{equation}
remains valid. Tackling this problem analytically turns out to be extremely challenging, see the proof of Theorem \ref{main_result} in \cite{Schneider}. \\

In our simulations, we consider sequences $(n_k, p_k)_k$ defined by $n_k = w\sqrt[\alpha]{k}$ and $p_k = \mu / n_k$ for parameters $w, \mu > 0$. The values of $w$ and $\mu$ should, ideally, not matter for the asymptotics and thus for the pursuit of $\alpha^*$. Suitable choices of $w$ and $\mu$ for given $\alpha$ are still necessary for practical reasons to ensure that the asymptotic behavior becomes visible for the values of $k$ covered by the simulations. 
For any selection $(\alpha, w, \mu)$, we calculate the posterior probability of the true parameter $n_k$ and the MSE of different estimators for values of $k$ up to $10^{11}$. In order to achieve these extremely large observation numbers, we take care to minimize the number of operations when expressing the beta-binomial likelihood $L_{a,b}$ in our implementation. Since $L_{a,b}$ does not depend on the order of the observations but only on the frequencies of each distinct outcome $x_i$, the runtime depends on $n_k$ (the number of different values that $x_i$ can take) instead of $k$ itself.

\begin{figure}
  \centering
  {\scriptsize\input{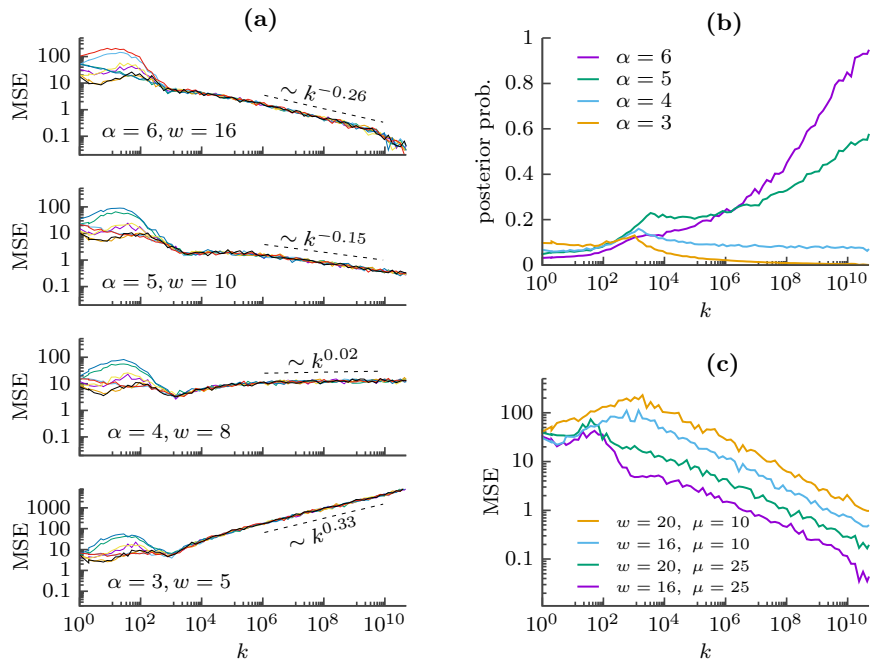}}
  \caption{Asymptotic behavior of the scale estimator and posterior contraction. \textbf{(a)} shows log-log plots of the MSE of several scale estimators in different asymptotic scenarios $(\alpha, w, \mu)$. The value $\mu$ was set to $25$ in each simulation, and the parameters for the scale estimators were picked as all possible combinations of $\gamma\in\{0.5, 1\}$, $a\in\{1, 5\}$, and $b\in\{1, 5\}$. \textbf{(b)} shows the empirical mean of the posterior probabilities $\Pi(N = n_k^0 \,|\, \mathbf{X}^k)$ for the same four settings depicted in (a). \textbf{(c)} shows the MSE of the scale estimator with parameters $\gamma = a = b = 1$ for constant $\alpha = 6$ and varying values of $w$ and $\mu$.}
  \label{fig:asymptotics}
\end{figure}

Figures \ref{fig:asymptotics}a--b show the (empirical) mean posterior probability in \eqref{eq:contraction} and the (empirical) mean square error (MSE) between $\hat{n}$ and $n$ for different scale estimators $\hat{n}$ in several scenarios $(\alpha, w, \mu)$. The number of samples was set to 200. It is clearly visible that the choice $\alpha = 6$ leads to posterior consistency (which is in good agreement with Theorem \ref{main_result}), since the posterior probability approaches 1 while the MSE converges to 0. However, the simulations indicate that this also holds true for $\alpha = 5$. For $\alpha = 4$, it becomes questionable whether posterior contraction will eventually happen. The choice $\alpha = 3$, in contrast, leads to a clear increase of the MSE with increasing $k$, and posterior contraction evidently fails.

\begin{figure}
  \centering
  {\scriptsize\input{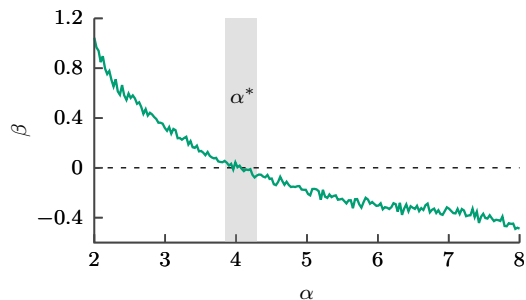}}
  \caption{Relation between $\alpha$ and $\beta$. For a given order $\alpha$, the corresponding value of $\beta$ was determined by conducting simulations like in Figure \ref{fig:asymptotics}a and fitting the slope for $k$ between $10^7$ and $10^9$. The graph shows that the zero point $\alpha^*$ of the conjectured function $\beta(\alpha)$ has to be in the vicinity of $4$.}
  \label{fig:slopes}
\end{figure}

An interesting observation is the power law behavior $\sim k^{-\beta}$ of the MSE, which is revealed by linear segments in the respective log-log plots. Figure \ref{fig:asymptotics}a shows that the slope $\beta$ is independent of the chosen estimator, and \ref{fig:asymptotics}c suggests that it might also be independent of $w$ and $\mu$. We can therefore consider $\beta$ as a function $\beta(\alpha)$ of $\alpha$ alone. A numerical approximation of $\alpha^*$ is then given by the value of $\alpha$ where $\beta$ changes sign, i.e.,
\begin{equation*}
  \beta(\alpha^*) = 0.
\end{equation*}
Since $\beta(\alpha)$ is strictly monotone, as a higher number $k$ of observations will lead to better estimates, such an $\alpha^*$ is uniquely defined. Figure \ref{fig:slopes} displays an approximation of the graph of $\beta(\alpha)$ for values between $\alpha = 2$ and $\alpha = 8$. The respective slopes are estimated by linear least squares regressions for $k$ between $10^7$ and $10^{9}$. Even though our numerical results do not allow us to establish the precise functional relation between $\alpha$ and $\beta$, it becomes clear that $\alpha^*$ indeed has to be close to $4$.\\

\begin{figure}
  \centering
  {\scriptsize\input{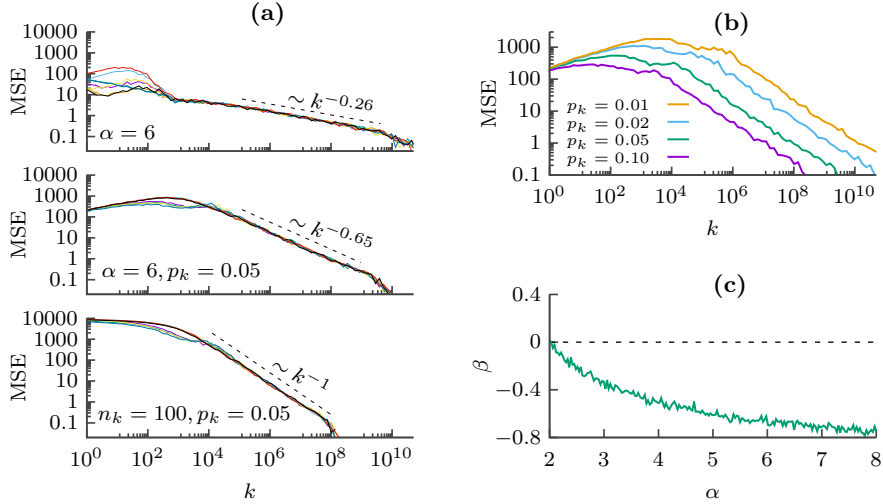}}
  \vspace{0.1cm}
  \caption{Comparison of alternative asymptotic settings. \textbf{(a)} shows the MSE for three different asymptotic scenarios. In the first plot, $n_k$ and $p_k$ behave like in Figure \ref{fig:asymptotics} with $w = 16$ and $\mu = 25$. In the second plot, $p_k$ is fixed to the value $0.05$, while $n_k$ still increases with $k$ ($w = 16$). The third plot addresses the scenario where both $n_k$ and $p_k$ are held fixed. \textbf{(b)} shows the scenario of growing $n_k$ (with $\alpha = 6$ and $w = 16$) and different fixed values $p_k$. The graph shows that the slope in the linear segment does not depend on $p_k$. \textbf{(c)} shows the relation between $\beta$ and $\alpha$ for the scenario with fixed $p_k$ and growing $n_k$. The values of the slopes $\beta$ are determined as described in Figure \ref{fig:slopes}, with adapted ranges for $k$.}
  \label{fig:additional}
\end{figure}

For comparison, we additionally conducted simulations that target other asymptotic regimes. First, we keep $p_k$ constant and let $n_k$ again increase with the sample size, $n_k = w\sqrt[\alpha]{k}$. In this scenario, a properly rescaled binomial random variable converges to a standard normal distribution. Our simulations confirm that estimation of $n$ is easier in this case: the MSE in Figure \ref{fig:additional}a decreases faster when $\alpha=6$ and $p_k=0.05$ is fixed compared to $\alpha=6$ and $p_k\rightarrow0$. Since the rate of convergence $\beta$
in this alternative setting seems to be independent of the specific choice of $p_k = \mathrm{const}$, see Figure \ref{fig:additional}b, we can again look at the smallest order $\alpha$ that still exhibits consistency. Indeed, Figure \ref{fig:additional}c reveals that the estimation of $n$ remains consistent over a larger range of values for $\alpha$ in this setting, approximately as long as $\alpha>2$ (compared to $\alpha > 4$ in the original setting).

The last asymptotic regime we consider is the classical one for parameter estimation, where $n_k = n$ and $p_k = p$ both stay constant as $k$ grows to infinity. Figure \ref{fig:additional}a covers this regime in the last plot. It affirms that estimating $n$ is easiest in this setting, and we obtain the expected rate $\sim k^{-1}$ for the convergence of the MSE towards zero.

\section{Discussion}
\label{sec:4}

Theorem \ref{main_result} (see \cite{Schneider}) shows posterior contraction under diverging parameters $n_k$ and $p_k$ as long as $(n_k,p_k)\in \mathcal{M}_{\lambda}$, which implies $n_k=O(\sqrt[6]{k/\log(k)})$. The aim of our simulation study in Section \ref{sec:3} was to explore the minimal rate $\sqrt[\alpha]{k}$ for $n_k$ such that posterior consistency remains valid. The difference in the permissible rates turns out to be rather small, since
our investigation suggests that $\alpha=5$ still allows for consistent estimation, whereas $\alpha=3$ clearly leads to inconsistency. Figure \ref{fig:slopes} shows that the true boundary $\alpha^*$ is likely close to 4, indicating that Theorem \ref{main_result} cannot be improved fundamentally.

Several aspects of our simulations and findings deserve further commentary. First, Figure \ref{fig:asymptotics}c reveals that the slope $\beta$ is not strongly affected by the parameters $w$ and $\mu$ in the settings that we tested. However, our numerical approach is not suitable to verify questions like this with a high degree of confidence. For example, our numerics become instable for values $k > 10^{11}$.

Secondly, we additionally conducted simulations for other estimators than the scale estimator \eqref{E:ScaleBayes} that are not shown in the article. For example, we tested various versions of the Bayesian estimator given in \cite{Draper-Guttman}. While their performance for $k \le 10^3$ varies quite much -- similar to the different estimators shown in Figure \ref{fig:asymptotics}a -- their asymptotic performance is exactly the same as for the scale estimator. Notably, the maximum likelihood estimator also exhibits the very same asymptotic behavior, even though it performs poorly in the regime of smaller $k$. 
The sample maximum, in contrast, shows a completely different behavior: the MSE diverges even for $n_k \sim \log(k)$. This illustrates the sharpness of the assumptions for Lemma 10 in \cite{Schneider}, which states that the sample maximum is consistent if $n_k\log(n_k)<c\log(k)$ for $c<1$.

Finally, we consistently observed a phase transition in all simulations when the MSE drops below a value of about $0.1$, where it changes its behavior and begins to decreases faster than $\sim k^\beta$. Indeed, it seems to decay exponentially from that point on. We conjecture that this happens due to the discreteness of $n$, which means that the MSE cannot measure small deviations $|\hat{n} - n| < 1$ from the real $n$ without dropping to zero. Rather, if the posterior contracts so much that we estimate $n$ correctly most of the time, the MSE essentially captures the probability that $\hat{n}$ lies outside of the interval $(n-1,n+1)$, and such probabilities usually decay exponentially fast. 
For applications, the rate of the MSE before the exponential decay is often much more interesting. 
One instructive example in this context is the sample maximum in the setting of fixed $n$ and $p$, for which we know from Section \ref{sec:1} that it converges exponentially fast. However, as argued above, this only takes place when the MSE is already very small, and simulations suggest that the rate of convergence is much slower if the MSE is larger than $0.1$. For instance, if $p = 0.2$ and $n = 25$, we find $\beta \approx -0.13$. Thus, even though the true asymptotic behavior of the sample maximum is exponential, the practically meaningful rate of convergence is considerably worse than the rate $k^{-1}$ of the Bayesian estimators.

\section*{Acknowledgements}
Support of the DFG RTG 2088 (B4) and DFG CRC 755 (A6) is gratefully acknowledged.

%
%
%

\end{document}